\documentclass[12pt]{amsproc}
\usepackage{amsmath,amssymb,amscd,euscript,oldgerm}

\begin{document}

\title[Countable Degree-1 Saturation]{Countable Degree-1 Saturation of Certain $C^*$-Algebras Which Are Coronas of Banach Algebras}
\author{Dan-Virgil Voiculescu}
\address{D.V. Voiculescu \\ Department of Mathematics \\ University of California at Berkeley \\ Berkeley, CA\ \ 94720-3840}
\thanks{Research supported in part by NSF Grant DMS 1301727.}
\date{}

\begin{abstract}
We study commutants modulo some normed ideal of 
$n$-tuples of operators which satisfy a certain approximate unit condition relative to the ideal. We obtain results about the quotient of these Banach algebras by their ideal of compact operators being $C^*$-algebras which have the countable degree-1 saturation property in the model-theory sense of I.~Farah and B.~Hart. We also obtain results about quasicentral approximate units, multipliers and duality.
\end{abstract}

\maketitle

\setcounter{section}{-1}
\section{Introduction}
\label{sec0}

The study of the commutant modulo the Hilbert--Schmidt class of a normal operator with rich spectrum (\cite{13}, \cite{2}) has shown that this Banach algebra together with its ideal of compact operators resembles in many ways the pair consisting of the algebra ${\EuScript B}({\EuScript H})$ of all operators on a Hilbert space ${\EuScript H}$ and the ideal ${\EuScript K}({\EuScript H})$ of compact operators and that the analog of the Calkin algebra is also a $C^*$-algebra. The purpose of this paper is to further develop this analogy. On one hand, we go beyond the case of a normal operator \cite{13} or of a commuting $n$-tuple of hermitian operators \cite{2} and deal with a general non-commuting $n$-tuple of operators and its commutant modulo a normed ideal which satisfies a certain quasicentral approximate unit condition relative to the $n$-tuple. The main result we obtain in this general setting is that countable degree-1 saturation, in the model theory sense of (\cite{6}), holds for our analog of the Calkin algebra, which is still a $C^*$-algebra in this general case. In what follows most of the time we will refer to countable degree-1 saturation simply as ``degree-1 saturation'', for the sake of brevity. This adds to the list of nice properties of these analogs of the Calkin algebra and also adds to the list of $C^*$-algebras satisfying degree-1 saturation (\cite{6}). We also obtain a few other results, like for instance the existence of quasicentral approximate units for the ideal of compact operators in the Banach algebra we consider, as well as generalizations of some of the multiplier and duality results in \cite{13}.

Perhaps, the results on the analogues of the Calkin algebra which we obtain, give hope that these algebras may be a good place to apply extensions of bi-variant $K$-theory beyond $C^*$-algebras (\cite{5}) and cyclic cohomology (\cite{4}).

Besides the introduction this paper has five sections.

Section~\ref{sec1} deals with preliminaries. Especially, in preparation for the later sections, we recall certain basic facts about normed ideals of compact operators (\cite{9}, \cite{12}) and about the invariant $k_y({\EuScript I})$ where ${\EuScript I}$ is a normed ideal and $\tau$ an $n$-tuple of operators, which we used in our work on normed ideal perturbations of Hilbert space operators (\cite{15}, \cite{14}, \cite{17}).

The main result of Section~\ref{sec2} is the existence of quasicentral approximate units for the compact ideal of the Banach algebras we study. The construction we use has some of the flavor of the tridiagonal construction we used in our original proof of the non-commutative Weyl--von~Neumann theorem \cite{16}, before the concept of quasicentral approximate units was abstracted (\cite{1}, \cite{11}). The fact that the analogue of the Calkin algebra is a $C^*$-algebra is also in this section.

Section~\ref{sec3} gives the countable degree-1 saturation for the analogue of the Calkin algebra. The proof is along similar lines to those of the proof for coronas of $C^*$-algebras of Farah and Hart (\cite{6}) with the added technical difficulties arising from Banach algebra norms which don't allow continuous functional calculus. On the other hand, we were helped by the fact that in the case of the Calkin algebra the main technical lemma and the glueing construction simplify and becomes reminiscent of the tridiagonal construction and the kind of approximately commuting partition of unity used to glue parts of operators in \cite{16}.

Section~\ref{sec4} deals with generalizations of multiplier and duality results from \cite{13} to the general setting. Here once appropriate assumptions are found, the proofs in \cite{13} generalize immediately.

Section~\ref{sec5} is a section of concluding remarks.

The author gratefully acknowledges the opportunity to learn about degree-1 saturation from attending the $C^*$-algebra meeting at Oberwolfach in August 2013 and the lecture of Farah Ilijas at the meeting and subsequent discussions with him.

\section{Preliminaries}
\label{sec1}

Throughout this paper the term {\em normed ideal} will be used as an abbreviation for {\em symmetrically normed ideal} (\cite{9}, \cite{12}) of compact operators on a separable infinite dimensional complex Hilbert space ${\EuScript H}$. This is an ideal $0 \ne {\EuScript I} \subset {\EuScript B}({\EuScript H})$ of the algebra of all bounded operators on ${\EuScript H}$ which is contained in ${\EuScript K}({\EuScript H})$ the ideal of compact operators and which is endowed with a certain norm $|~|_{{\EuScript I}}$ with respect to which is a Banach space. The norm is given by $|T|_{{\EuScript I}} = |T|_{\Phi} = \Phi(s_1(T),s_2(T),\dots)$ where $\Phi$ is a norming function (see \S\ref{sec3} in \cite{9}) and $s_1(T) \ge s_2(T) \ge \dots$ are the $s$-numbers of $T$. Given a norming function $\Phi$ we will use the notation in \cite{9} and denote by $({\textgoth S}_{\Phi},|~|_{\Phi})$ and $({\textgoth S}_{\Phi}^{(0)},|~|_{\Phi})$ the normed ideals which are the set of all compact operators $T$ so that $|T|_{\Phi} < \infty$ and, respectively, the closure in ${\textgoth S}_{\Phi}$ of ${\EuScript R}({\EuScript H})$ the ideal of finite rank operators. We will always leave out ${\EuScript K}({\EuScript H})$ as a normed ideal in our considerations. If $({\EuScript I},|~|_{{\EuScript I}})$ is a normed ideal we shall also use the notation ${\EuScript I}^{(0)}$ for the closure of ${\EuScript R}({\EuScript H})$ in ${\EuScript I}$. Remark that since $|~|_{{\EuScript I}} = |~|_{\Phi}$ for some norming function $\Phi$, ${\EuScript I}^{(0)} = {\textgoth S}_{\Phi}^{(0)}$. Also if $|~|_{{\EuScript I}} = |~|_{\Phi}$ we clearly have ${\textgoth S}_{\Phi}^{(0)} \subset {\EuScript I} \subset {\textgoth S}_{\Phi}$ and if ${\textgoth S}_{\Phi}^{(0)} = {\textgoth S}_{\Phi}$ the function $\Phi$ is called ``mononorming'' \cite{9}.

If $\tau = (T_j)_{1 \le j \le n}$ is an $n$-tuple of operators the definition of the number
\[
k_{{\EuScript I}}(\tau) = \liminf_{A \in {\EuScript R}_1^+({\EuScript H})}|[A,\tau]|_{{\EuScript I}}
\]
from (\cite{15} see also \cite{14}, \cite{17}), where $({\EuScript I},|~|_{{\EuScript I}})$ is a normed ideal and ${\EuScript R}_1^+({\EuScript H}) = \{A \in {\EuScript R}({\EuScript H}) \mid 0 \le A \le I\}$ the lim~inf being with respect to the natural order on ${\EuScript R}_1^+({\EuScript H})$ and where we use the notation $[A,\tau] = ([A,T_j])_{1 \le j \le n}$ and $|(X_j)_{1 \le j \le n}|_{{\EuScript I}} = \max_{1 \le j \le n} |X_j|_{{\EuScript I}}$. If $|~|_{{\EuScript I}} = |~|_{\Phi}$ we also write $k_{\Phi}(\tau)$ for $k_{{\EuScript I}}(\tau)$.

We will be mainly interested here in the condition $k_{{\EuScript I}}(\tau) = 0$. Results concerning this are summarized in \cite{17}. For instance, if $\tau$ is an $n$-tuple of commuting hermitian operators and ${\EuScript I} = {\EuScript C}_n$ the Schatten von~Neumann class, then we have $k_{{\EuScript C}_n}(\tau) = 0$ if $n \ge 2$. This implies the fact that $k_{{\EuScript C}_2}(N) = 0$ if $N$ is a normal operator which underlies the results in \cite{13}.

We should also recall (see \cite{15} or \cite{14}) that $k_{{\EuScript I}}(\tau) = 0$ is equivalent to $k_{{\EuScript I}}(\tau\coprod \tau^*) = 0$ where $\tau^* = (T_j^*)_{1 \le j \le n}$ or to $k_{{\EuScript I}}(\text{Re } \tau \coprod \text{Im } \tau) = 0$ where $\text{Re } \tau = (\text{Re } T_j)_{1 \le j \le n}$ and $\text{Im } \tau = (\text{Im } T_j)_{1 \le j \le n}$. The condition $k_{{\EuScript I}}(\tau) = 0$ is also equivalent to the existence of a sequence $A_n \in {\EuScript R}_1^+({\EuScript H})$ such that $A_n \overset{w}{\longrightarrow} I$ and $|[A_n,\tau]|_{{\EuScript I}} \to 0$ as $n \to \infty$ or also to the existence of a sequence $A_n \uparrow I$, $A_n \in {\EuScript R}_1^+({\EuScript H})$ satisfying additional conditions like $m > n \Rightarrow A_mA_n = A_n$ and $A_nB_n = B_n$ where $B_n \in {\EuScript R}({\EuScript H})$ are given and so that $|[A_n,\tau]|_{{\EuScript I}} \to 0$ as $n \to \infty$.

\section{Approximate units}
\label{sec2}

Let $\tau = (T_j)_{1 \le j \le n}$, $T_j = T_j^*$, $1 \le j \le n$ be an $n$-tuple of hermitian operators in ${\EuScript B}({\EuScript H})$ and let $({\EuScript I},|~|_{{\EuScript I}})$ be a normed ideal we define
\[
{\EuScript E}(\tau;{\EuScript I}) = \{X \in {\EuScript B}({\EuScript H}) \mid [X,T_j] \in {\EuScript I},\ 1 \le j \le n\}
\]
and ${\EuScript K}(\tau;{\EuScript I}) = {\EuScript E}(\tau;{\EuScript I}) \cap {\EuScript K}({\EuScript H})$. Then ${\EuScript E}(\tau;{\EuScript I})$ is a Banach algebra with the norm $\||X\|| = \|X\| + |[X,\tau]|_{{\EuScript I}}$ with an isometric involution $\||X^*\|| = \||X\||$ and ${\EuScript K}(\tau;{\EuScript I})$ is a closed two-sided ideal, which is also closed under the involution. We shall denote by ${\EuScript P}({\EuScript H})$ the finite-rank hermitian projections. Clearly ${\EuScript P}({\EuScript H}) \subset {\EuScript R}({\EuScript H}) \subset {\EuScript K}(\tau;{\EuScript I})$.

\bigskip
\noindent
{\bf 2.1. Proposition.} {\em Assume $k_{{\EuScript I}}(\tau) = 0$.

{\em a)} If $P \in {\EuScript P}({\EuScript H})$ and $\epsilon > 0$, then there is $A \in {\EuScript R}_1^+({\EuScript H})$ so that $P \le A$ and $\||A\|| < 1 + \epsilon$.

{\em b)} If ${\EuScript R}({\EuScript H})$ is dense in ${\EuScript I}$ and $P \in {\EuScript P}({\EuScript H})$, $K_r \in {\EuScript K}(\tau;{\EuScript I})$, $1 \le r \le m$ and $\epsilon > 0$, then there is $A \in {\EuScript R}_1^+({\EuScript H})$ so that $P \le A$, $\||(I-A)K_r\|| < \epsilon$, $1 \le r \le m$ and $\||A\|| < 1 + \epsilon$.
}

\bigskip
\noindent
{\bf {\em Proof.}} a) Since $k_{{\EuScript I}}(\tau) = 0$ there is $A \in {\EuScript R}_1^+({\EuScript H})$, $P \le A$ so that $|[A,\tau]|_{{\EuScript I}} < \epsilon$, which in view of the fact that $\|A\| \le 1$ gives $\||A\|| < 1 + \epsilon$.

b) Since $[K_r,T_j] \in {\EuScript I}$, $1 \le r \le m$, $1 \le j \le n$ and ${\EuScript R}({\EuScript H})$ is dense in ${\EuScript I}$, there is a projection $Q \in {\EuScript P}({\EuScript H})$ so that $|(I-Q)[K_r,T_j]|_{{\EuScript I}} < \epsilon/4$ and $\|(I-Q)K_r\| < \epsilon/4$, $1 \le r \le m$, $1 \le j \le n$. Clearly, we may assume without loss of generality that $P \ge Q$ and $\||K_r\|| \le 1$, $1 \le r \le m$. Using a), there is $A \in {\EuScript R}_1^+({\EuScript H})$ so that $Q \le P \le A$ and $|[A,\tau]|_{{\EuScript I}} < \epsilon/4$. We have
\[
\|(I-A)K_r\| \le \|(I-Q)K_r\| < \epsilon/4
\]
and
\[
\begin{aligned}
|[(I-A)K_r,\tau]|_{{\EuScript I}} &\le |[A,\tau]|_{{\EuScript I}}\|K_r\| + \max_{1 \le j \le n} |(I-A)[K_r,T_j]|_{{\EuScript I}} \\
&< \epsilon/4 + \epsilon/4 = \epsilon/2.
\end{aligned}
\]
It follows that $\||(I-A)K_r\|| < \epsilon$.\qed

\bigskip
\noindent
{\bf 2.2. Corollary.} {\em If $K_{{\EuScript I}}(\tau) = 0$ and ${\EuScript R}({\EuScript H})$ is dense in ${\EuScript I}$, then ${\EuScript R}({\EuScript H})$ is dense in ${\EuScript K}(\tau;{\EuScript I})$.
}

\bigskip
\noindent
{\bf 2.3. Proposition.} {\em Assume $k_{{\EuScript I}}(\tau) = 0$ and ${\EuScript I}^{(0)}={\EuScript I}$, that is ${\EuScript R}({\EuScript H})$ is dense in ${\EuScript I}$. Let $X_1,\dots,X_m \in {\EuScript E}(\tau;{\EuScript I})$, $K_1,\dots,K_r \in {\EuScript K}(\tau;{\EuScript I})$, $P \in {\EuScript P}({\EuScript H})$ and $\epsilon  > 0$ be given. Then there is $B \in {\EuScript R}_1^+({\EuScript H})$ so that $P \le B$, $\||B\|| < 1 + \epsilon$
\[
\||(I-B)K_j\|| < \epsilon,\ \||[X_p,B]\|| < \epsilon
\]
for $1 \le j \le r$, $1 \le p \le m$.
}

\bigskip
\noindent
{\bf {\em Proof.}} Without loss of generality we will assume that $X_p = X_p^*$, $1 \le p \le m$. Since ${\EuScript I} = {\EuScript I}^{(0)}$ there is $P_0 \in {\EuScript P}({\EuScript H})$ so that $P \le P_0$ and
\[
|(I-P_0)[X_p,\tau]|_{{\EuScript I}} + |[X_p,\tau](I-P_0)|_{{\EuScript I}} < \epsilon/2.
\]
Applying repeatedly Proposition~2.1 we can find $P_s \in {\EuScript P}({\EuScript H})$, $A_s \in {\EuScript R}_1^+({\EuScript H})$,
\[
\begin{aligned}
P_0 &\le P_1 \le P_2 \le \dots, \\
A_0 &\le A_1 \le A_2 \le \dots
\end{aligned}
\]
so that $P_s \uparrow I$ as $s \to \infty$ and $P_s \le A_s \le P_{s+1}$, $(I-P_{s+1})X_pA_s = 0$, $(I-P_{s+1})T_lA_s = 0$, $(I-P_{s+1})X_pA_s = 0$ (that is $P_{s+1}{\EuScript H} \supset X_pA_s{\EuScript H} + T_lA_s{\EuScript H}$), $\||A_s\|| < 1 + \epsilon 2^{-s-1}$ and $\||(I-A_s)K_j\|| < \epsilon$ for $1 \le p \le m$, $1 \le l \le n$, $1 \le j \le r$ and all $s \ge 0$. Let $B = N^{-1}(A_1+\dots + A_N)$. We will show that choosing $N$ sufficiently large, $B$ will have all the desired properties. Clearly, since $A_s \ge P$, $1 \le s \le N$ we will also have the same inequality for their mean, that is $B \ge P$. Similarly, $(I-B)K_j$ is the mean of the $(I-A_s)K_j$, $1 \le s \le N$ and this gives $\||(I-B)K_j\|| < \epsilon$. Also, the same kind of argument gives $\||B\|| < 1 + N^{-1}\epsilon$.

To prove that $\||[X_p,B]\|| < \epsilon$ if $N$ is large enough we will show that $\|[X_p,B]\| \to 0$ and $|[X_p,B],\tau]|_{{\EuScript I}} \to 0$ as $N \to \infty$. Remark that the conditions on $P_s$, $A_s$, $X_p$, $T_l$ imply that in the orthogonal sum decomposition
\[
{\EuScript H} = P_0{\EuScript H} \oplus (P_1-P_0){\EuScript H} \oplus (P_2-P_1){\EuScript H} \oplus \dots
\]
we have that $A_s$ is block-diagonal, while the $X_p$ and $T_l$, being hermitian, are block-tridiagonal. With the notation $Q_0 = P_0$, $Q_s = P_s - P_{s-1}$, $s \ge 1$, we have $A_{s-1} = Q_0 + \dots + Q_{s-1} + Q_sA_{s-1}Q_s$ if $s \ge 1$. It follows that
\[
\begin{split}
\left\|B-\left(Q_0 + \sum_{1 \le s \le N}\left( 1 - \frac {s-1}{N}\right) Q_s\right)\right\|  \\
= \left\| N^{-1} \sum_{1 \le s \le N} Q_{s+1}A_sQ_{s+1}\right\| 
\le N^{-1}.
\end{split}
\]
Hence the tridiagonality gives
\[
\begin{aligned}
\|[B,X_p]\| &\le 2N^{-1}\|X_p\| + \left\|\left[ Q_0 + \sum_{1 \le s \le N} \left( 1 - \frac {s-1}{N}\right)Q_s,X_p\right]\right\| \\
&\le 2N^{-1} + N^{-1}\left( \left\| \sum_{1 \le s \le N} Q_{s+1}X_pQ_s\right\| + \left\| \sum_{1 \le s \le N} Q_sX_pQ_{s+1}\right\|\right) \\
&\le 4N^{-1}\|X_p\|
\end{aligned}
\]
and hence $\|[B,X_p]\| \to 0$ as $N \to \infty$.

Since we may choose $P_0 \ne 0$, we have $\|B\| = 1$ and hence $\||B\|| < 1 + N^{-1}\epsilon$ gives $|[B,\tau]|_{{\EuScript I}} < \epsilon N^{-1}$. It follows that
\[
\begin{aligned}
|[[B,X_p],\epsilon]|_{{\EuScript I}} &\le 2|[B,\tau]|_{{\EuScript I}}\|X_p\| + |[B,[X_p,\tau]]|_{{\EuScript I}} \\
&\le 2N^{-1}\epsilon\|X_p\| + |(I-B)[X_p,\tau]|_{{\EuScript I}} + |[X_p,\tau](I-B)|_{{\EuScript I}}.
\end{aligned}
\]
Since $B \ge P_0$, it follows that
\[
|(I-B)[X_p,\tau]|_{{\EuScript I}} + |[X_p,\tau](I-B)|_{{\EuScript I}} < \epsilon/2.
\]
Hence $|[B,X_p],\tau|_{{\EuScript I}} < \epsilon$ for $N$ large enough.\qed

\bigskip
Remark that since the assumptions ${\EuScript I}^{(0)} = {\EuScript I}$, $k_{{\EuScript I}}(\tau) = 0$ imply that ${\EuScript R}({\EuScript H})$ is dense in ${\EuScript K}(\tau;{\EuScript I})$ is a separable Banach space and applying repeatedly Proposition~2.3 we immediately can tie it a somewhat stronger form, which we state as the next corollary.

\bigskip
\noindent
{\bf 2.4. Corollary.} {\em Assume $k_{{\EuScript I}}(\tau) = 0$ and ${\EuScript I}^{(0)} = {\EuScript I}$. Let $X_1,\dots,X_m \in {\EuScript E}(\tau;{\EuScript I})$ and a sequence $Y_s \in {\EuScript R}({\EuScript H})$, $s \in {\mathbb N}$ be given. Then there is a sequence $A_s \in {\EuScript R}_1^+({\EuScript H})$ so that $A_sY_s = Y_s$ and $A_sA_t = A_t$, $A_sX_pA_t = X_pA_t$ if $s > t$ and moreover
\[
A_s \uparrow I,\ \||A_s\|| \to 1,\ \||(I-A_s)K\|| \to 0,\ \||[X_p,A_s]\|| \to 0
\]
as $s \to \infty$ for all $K \in {\EuScript K}(\tau;{\EuScript I})$ and $1 \le p \le m$.
}

\bigskip
We pass now to the quotient Banach algebra with involution \linebreak ${\EuScript E}(\tau;{\EuScript I})/{\EuScript K}(\tau;{\EuScript I})$ which we shall denote by ${\EuScript E}/{\EuScript K}(\tau;{\EuScript I})$. If
\[
p: {\EuScript B}({\EuScript H}) \to {\EuScript B}({\EuScript H})/{\EuScript K}({\EuScript H}) = {\EuScript B}/{\EuScript K}({\EuScript H})
\]
is the canonical homomorphism to the Calkin algebra, which we shall denote by ${\EuScript B}/{\EuScript K}({\EuScript H})$, then there is a canonical isomorphism of ${\EuScript E}/{\EuScript K}(\tau;{\EuScript I})$ and the sub algebra $p({\EuScript E}(\tau;{\EuScript I}))$ of ${\EuScript B}/{\EuScript K}(\tau;{\EuScript I})$. In view of this we shall often also denote by $p$ the homomorphism ${\EuScript E}(\tau;{\EuScript I}) \to {\EuScript E}/{\EuScript K}(\tau;{\EuScript I})$, an admissible abuse of notation.

\bigskip
\noindent
{\bf 2.5. Proposition.} {\em We assume $k_{{\EuScript I}}(\tau) = 0$ and ${\EuScript I}^{(0)} = {\EuScript I}$. Given $X \in {\EuScript E}(\tau;{\EuScript I})$ and $\epsilon > 0$ there is $A \in {\EuScript R}_1^+({\EuScript H})$ so that $\|A\| = 1$, $\||A\|| < 1 + \epsilon$ and $\||(I-A)X\|| < \|p(X)\| + \epsilon$ where the norm of $p(X)$ is the ${\EuScript B}/{\EuScript K}({\EuScript H})$ norm. In particular, the norm of $X = {\EuScript K}(\tau;{\EuScript I})$ in ${\EuScript E}/{\EuScript K}(\tau;{\EuScript I})$ equals the norm of $p(X)$ in ${\EuScript B}/{\EuScript K}({\EuScript H})$. Thus algebraic embedding of ${\EuScript E}/{\EuScript K}(\tau;{\EuScript I})$ into ${\EuScript B}/{\EuScript K}({\EuScript H})$ is isometric and ${\EuScript E}/{\EuScript K}(\tau;{\EuScript I})$ identifies with a $C^*$-subalgebra of ${\EuScript B}/{\EuScript K}({\EuScript H})$.
}

\bigskip
\noindent
{\bf {\em Proof.}} We have stated this fact which is an immediate generalization of results in \cite{13} and \cite{3}, with a lot of detail, since it will be often used in the rest of this paper.

In view of our assumption, that $k_{{\EuScript I}}(\tau) = 0$, there are $A_n \uparrow I$, $A_n \in {\EuScript R}_1^+({\EuScript H})$ so that $|[A_n,\tau]|_{{\EuScript I}} \to 0$ as $n \to \infty$. Then also $\|(I-A_n)X\| \to \|p(X)\|$ as $n \to \infty$. We also have
\[
|[(I-A_n)X,\tau]|_{{\EuScript I}} \le |[A_n,\tau]|_{{\EuScript I}} \|X\| + |(I-A_n)[X,\tau]|_{{\EuScript I}}
\]
and the first term in the right-hand side $\to 0$ as $n \to \infty$ by the properties of the $A_n$, while the second also $\to 0$ since ${\EuScript I}^{(0)} = {\EuScript I}$ and $(I-A_n)[X,\tau]$ converges weakly to $0$ as $n \to \infty$.

The rest of the statement is well explained in the statement of the corollary itself.\qed

\section{Countable degree-1 saturation}
\label{sec3}

In this section we prove what amounts to countable degree-1 saturation of ${\EuScript E}/{\EuScript K}(\tau;{\EuScript I})$ under the assumption that $k_{{\EuScript I}}(\tau) = 0$, in the model-theory terminology of \cite{3}, \cite{6}. The result is given in Theorem~3.3, which is formulated in operator-algebra terms, using one of the equivalent definitions of countable degree-1 saturation which can be found in \cite{6}.

We also remind the reader that the adjective ``countable'' will be omitted most of the time.

We begin with a rather standard technical fact, which we record as the next lemma.

\bigskip
\noindent
{\bf 3.1. Lemma.} {\em Let $G = G^* \in {\EuScript E}(\tau;{\EuScript I})$ be such that $\||G - \frac {3}{2} I\|| \le 1$. Then $G^{1/2} \in {\EuScript E}(\tau;{\EuScript I})$ and there is a universal constant $C$, so that
$\||[G^{1/2},X]\|| \le C\||[G,X]\||$ if $X \in {\EuScript E}(\tau;{\EuScript I})$ and $|[G^{1/2},\tau]|_{{\EuScript I}} \le C|[G,\tau]|_{{\EuScript I}}$.
}

\bigskip
\noindent
{\bf {\em Proof.}} The Lemma is an easy consequence of the functional calculus formula
\[
G^{1/2} = (2\pi i)^{-1} \int_{\Gamma} (zI-G)^{-1}z^{1/2}dz,
\]
where $\Gamma$ is the circle $|z-3/2| = 5/4$, and of the fact that for $z \in \Gamma$ we have:
\[
\begin{aligned}
\||(zI-G)^{-1}\|| &= 4/5\||(4/5(z-3/2)I - 4/5(G-3/2I))^{-1}\|| \\
&\le (1-4/5)^{-1} = 5
\end{aligned}
\]
and
\[
[(zI-G)^{-1},X] = (zI-G)^{-1}[G,X](zI-G)^{-1}.
\]
\qed

\bigskip
\noindent
{\bf 3.2. Lemma.} {\em Assume $k_{{\EuScript I}}(\tau) = 0$ and ${\EuScript I} = {\EuScript I}^{(0)}$. Let $M_n \in {\EuScript E}(\tau;{\EuScript I})$, $n \in {\mathbb N}$, $\epsilon_m \downarrow 0$ as $m \to \infty$, $P_k \in {\EuScript P}({\EuScript H})$, $P_k \uparrow I$ as $k \to \infty$ and an increasing function $\varphi: {\mathbb N} \to {\mathbb N}$ be given. Then there are $R_m \in {\EuScript R}_1^+({\EuScript H})$, $m \in {\mathbb N}$ so that

$1^0$. $\displaystyle{\sum_{m \ge 1} R_m^2 = I}$

$2^0$. the $R_m$'s commute

$3^0$. $\|R_m\| = 1$ and $|[R_m,\tau]|_{{\EuScript I}} < \epsilon_m$ if $m \ge 2$

$4^0$. $R_mP_n = 0$ if $m \ge n + 2$, $n \ge 1$

$5^0$. $\||[R_m,M_k]\|| < \epsilon_m$ if $k \le \varphi(m)$, $m \ge 2$

$6^0$. $R_nM_kR_m = 0$, $R_nR_m = 0$

\noindent
if $k \le \varphi(m)$, $k \le \varphi(n)$, $|n-m| \ge 2$, $m \ge 2$, $n \ge 2$.
}

\bigskip
\noindent
{\bf {\em Proof.}} There will be no loss of generality to assume that $M_1 = I$ and $M_k = M_k^*$, $k \in {\mathbb N}$. Given $\delta_m < 1/10$ $\delta_m \downarrow 0$ $m \to \infty$, we can use Proposition~2.1 repeatedly, to find a sequence of projections $E_k \in {\EuScript P}({\EuScript H})$, $E_k \uparrow I$ as $k \to \infty$ and a sequence of $A_k \in {\EuScript R}_1^+({\EuScript H})$, $A_k \uparrow I$ as $k \to \infty$, $A_1 = 0$ satisfying the following conditions:
\[
E_k \ge P_k,\ E_k \ge A_k,\ (I-E_k)M_pA_k = 0 \text{ if } p \le \varphi(k+2)
\]
and also
\[
\begin{aligned}
&A_{k+1} \ge E_k,\ \|A_{k+1}-E_k\| = 1, \\
&\||A_{k+1}\|| < 1 + \delta_{k+1} \text{ and} \\
&\||[A_{k+1},M_p]\|| < \delta_{k+1} \text{ if } p \le \varphi(k+2).
\end{aligned}
\]
Note that since the $E_k$'s are projections, $0 = A_1 \le E_1 \le A_2 \le E_2 \le \dots \le E_k \le A_{k+1} \le E_{k+1} \le \dots \le I$ implies that $\{A_k \mid k \ge 1\} \cup \{E_k \mid k \ge 1\}$ is a set of commuting operators.

With these preparations one would be tempted to define $R_n$ to be $(A_n-A_{n-1})^{1/2}$, however this would lead to difficulties with commutators and $\||\cdot\||$-norms because the square-root function is not differentiable at $0$. The remedy is to replace the $A_n$'s with another sequence $B_n$, $n \ge 1$, with roughly the same properties and which additionally is so designed that the problems with the square-roots are circumvented. We define
\[
B_n = I - (I-A_n^2)^2 = A_n^2(2I-A_n^2).
\]
Then $E_n \le A_{n+1} \le E_{n+1}$ easily gives $E_n \le B_{n+1} \le E_{n+1}$ and $A_1 = 0$ gives $B_1 = 0$. It is also easily seen that defining
\[
\begin{aligned}
R_n &= (B_n - B_{n-1})^{1/2},\ n \ge 2 \\
R_1 &= 0
\end{aligned}
\]
we have
\[
\begin{aligned}
B_n^{1/2} &= A_n(2I - A_n^2)^{1/2} \\
(I-B_n)^{1/2} &= I - A_n^2 \\
R_n &= (B_n(I-B_{n-1}))^{1/2} = B_n^{1/2}(I-B_{n-1})^{1/2} \\
&= A_n(2I-A_n^2)^{1/2}(I-A_{n-1}^2).
\end{aligned}
\]
Then for $n \ge 2$ we have
\[
\begin{aligned}
\||(2I-A_n^2)-3/2I\|| &= \||A_n^2-1/2I\|| \\
&= \|A^2-1/2I\| + |[A_n^2,\tau]|_{{\EuScript I}} \\
&\le 1/2 + 2\|A_n\| |[A_n,\tau]|_{{\EuScript I}} \le 1/2 + 2\delta_n \le 1.
\end{aligned}
\]
We can then apply Lemma~3.1 to $G = 2I - A_n^2$ and $X \in {\EuScript E}(\tau;{\EuScript I})$ and get that
\[
|[(2I-A_n^2)^{1/2},\tau]|_{{\EuScript I}} \le C|[A_n^2,\tau]|_{{\EuScript I}} \le 2C|[A_n,\tau]|_{{\EuScript I}} \le 2C\delta_n
\]
and
\[
\begin{aligned}
\||[(2I-A_n^2)^{1/2},X]\|| &\le C\||[A_n^2,X]\|| \\
&\le 2C\||[A_n,X]\||(1+\delta_n) \le 3C\||[A_n,X]\||.
\end{aligned}
\]

Remark that since $P_m \uparrow I$ we have $E_m \uparrow I$ and hence $B_m \uparrow I$ as $m \to \infty$. In view of $B_1 = 0$ we get that condition $1^0$ is satisfied by the $R_m$. Also, since the $B_m$ commute, it follows that the $R_m$ commute and that condition $2^0$ holds. Further, since $\|A_{k+1}-E_k\| = 1$, we have that $A_{k+1}$ has an eigenvector for the eigenvalue $1$ in $(E_{k+1}-E_k){\EuScript H}$, which is then also an eigenvector for the eigenvalue $1$ of $B_{k+1}$ and an eigenvector for the eigenvalue $0$ of $B_k$, so that it is an eigenvector for the eigenvalue $1$ for $(B_{k+1}-B_k)^{1/2} = R_{k+1}$. Thus we have $\|R_m\| = 1$ if $m \ge 2$, which is the first part of condition $2^0$. Further $R_mE_{m-2} = 0$ gives $R_mP_{m-2} = 0$, so that $4^0$ is satisfied. Also, since $M_1 = I$ and $M_k = M_k^*$, to check that $6^0$ holds, it suffices to check that $R_nM_kR_m = 0$ if $n \ge m+2$, $1 \le k \le \varphi(m)$. Indeed, we have $R_nM_kR_m = B_n^{1/2}(I-B_{n-1})^{1/2}M_kB_m^{1/2}(I-B_{m-1})^{1/2}$ and thus it suffices to show that $(I-B_{n-1})M_kB_m = 0$ if $1 \le k \le \varphi(m)$, $n \ge m+2$, which in turn will follow if we show that $(I-A_{n-1})M_kA_m = 0$ if $1 \le k \le \varphi(m)$, $n \ge m+2$. Note further that if $n \ge m+2$ we have $A_{n-1} \ge E_{n-2} \ge E_m$ and it suffices if $(I-E_m)M_kA_m$ for $k \le \varphi(m)$ for $m \ge 2$, which is satisfied in view of the construction of the $E_m$ and $A_m$.

We are thus left with having to deal with the second part of $3^0$ and $5^0$.

We have
\[
\begin{aligned}
|[R_m\tau]|_{{\EuScript I}} &= |[A_m(2I-A_m^2)^{1/2}(I-A_{m-1}^2),\tau]|_{{\EuScript I}} \\
&\le 2|[A_m,\tau]|_{{\EuScript I}} + |[(2I-A_m^2)^{1/2},\tau]|_{{\EuScript I}} + 2|[A_{m-1}^2,\tau]|_{{\EuScript I}} \\
&\le 2\delta_m + 2C\delta_m + 4\delta_{m-1}.
\end{aligned}
\]
Hence, choosing the $\delta_m$'s so that $2(C+1)\delta_m + 4\delta_{m-1} < \epsilon_m$ will insure that the second part of $3^0$ holds.

Turning to condition $5^0$, we have
\[
\begin{aligned}
\||[R_m,M_k]\|| &= \||[A_m(2I-A_m^2)^{1/2}(I-A_{m-1}^2),M_k]\|| \\
&\le \||[A_m,M_k]\||\,\||(2I-A_m^2)^{1/2}\||\,\||I-A_{m-1}^2\|| \\
&+ \||[(2I-A_m^2)^{1/2},M_k]\||\,\||A_m\||\,\||I-A_{m-1}^2\|| \\
&+ 2\||[A_{m-1},M_k]\||\,\||A_{m-1}\||\,\||A_m\||\,\||(2I-A_m^2)^{1/2}\|| \\
&\le \delta_m(2+|[(2I-A_m^2)^{1/2},\tau]|_{{\EuScript I}})(2+\delta_{m-1})^2 \\
&+ 3C\||[A_m,M_k]\||(1+\delta_m)(2+\delta_{m-1})^2 \\
&+ 2\delta_{m-1}(1+\delta_{m-1})(1+\delta_m)(2+|[(2I-A_m^2)^{1/2},\tau]|_{{\EuScript I}}) \\
&\le \delta_m(2+2C\delta_m)(2+\delta_{m-1})^2 \\
&+ 3C\delta_m(1+\delta_m)(2+\delta_{m-1})^2 \\
&+ 2\delta_{m-1}(1+\delta_{m-1})(1+\delta_m)(2+2C\delta_m),
\end{aligned}
\]
if $m \ge 2$ and $k \le \varphi(m)$. Clearly, the $\delta_m$'s can be chosen so that $\||[R_m,M_k]\|| < \epsilon_m$.\qed

\bigskip
By Proposition~2.5, ${\EuScript E}/{\EuScript K}(\tau;{\EuScript I})$ under the assumptions that $k_{{\EuScript I}}(\tau) = 0$ and ${\EuScript I} = {\EuScript I}^{(0)}$ is a $C^*$-algebra, actually a $C^*$-subalgebra of the Calkin algebra. Recall also that $p$ will denote both the homomorphism ${\EuScript B}({\EuScript H}) \to {\EuScript B}/{\EuScript K}({\EuScript H})$ as well as the homomorphism ${\EuScript E}(\tau;{\EuScript I}) \to {\EuScript E}/{\EuScript K}(\tau;{\EuScript I})$, which can be viewed as its restriction to ${\EuScript E}(\tau;{\EuScript I})$ (see the discussion preceding Proposition~2.3 and Proposition~2.5).

Let $X_j$, $X_j^*$, $j \in {\mathbb N}$ be non-commuting indeterminates and let
\[
f_n(X_1,\dots,X_n) = e_n + \sum_{j=1}^n(a_{jn}X_jb_{jn}+c_{jn}X_j^*d_{jn})
\]
where $e_n,a_{1n},\dots,a_{nn},\dots,b_{1n},\dots,b_{nn}$, $c_{1n},\dots,c_{nn},d_{1n},\dots,d_{nn}$ are in ${\EuScript E}/{\EuScript K}(\tau;{\EuScript I})$ so that the $f_n$ are non-commutative polynomials with coefficients not commuting with the variables. We shall denote the ring of such polynomials by ${\EuScript E}/{\EuScript K}(\tau;{\EuScript I})\langle X_j,X_j^* \mid j \in {\mathbb N}\rangle$, the $f_n$'s being polynomials of degree $\le 1$ in the indeterminates.

\bigskip
\noindent
{\bf 3.3. Theorem.} {\em Assume $k_{{\EuScript I}}(\tau) = 0$ and ${\EuScript I} = {\EuScript I}^{(0)}$. Let $e_n,a_{jn},b_{jn},c_{jn},d_{jn} \in {\EuScript E}/{\EuScript K}(\tau;{\EuScript I})$, $1 \le j \le n$, $n \in {\mathbb N}$, be such that there are $y_{jn} \in {\EuScript E}/{\EuScript K}(\tau;{\EuScript I})$, $1 \le j \le n$, $n \in {\mathbb N}$, so that $\|y_{jn}\| < 1$ and $|\|e_n + \sum_{1 \le j \le n} (a_{jn}y_{jm}b_{jn}+c_{jn}y_{jm}^*d_{jn})\| - r_n| < 1/m$ if $1 \le n \le m$, where $r_n \in {\mathbb R}$. Then there are $y_j \in {\EuScript E}/{\EuScript K}(\tau;{\EuScript I})$, $j \in {\mathbb N}$, so that $\|y_j\| \le 1$, for $j \in {\mathbb N}$ and
\[
\left\|e_n + \sum_{1 \le j \le n}(a_{jn}y_jb_{jn} + c_{jn}y_j^*d_{jn})\right\| = r_n
\]
for all $n \in {\mathbb N}$.
}

\bigskip
\noindent
{\bf {\em Proof.}} Let $E_n,A_{jn},B_{jn},C_{jn},D_{jn},Y_{jn} \in {\EuScript E}(\tau;{\EuScript I})$ for $1 \le j \le n$, $n \in {\mathbb N}$, be so that $p(E_n) = e_n$, $p(A_{jn}) = a_{jn}$, $p(B_{jn}) = b_{jn}$, $p(C_{jn}) = c_{jn}$, $p(D_{jn}) = d_{jn}$, $p(Y_{jn}) = y_{jn}$, $\||Y_{jn}\|| < 1$ and $|[Y_{jn},\tau]|_{{\EuScript I}} < \epsilon_n$ for some given sequence $\epsilon_n \downarrow 0$. It will be convenient to also introduce $f_n(X_1,\dots,X_n) \in {\EuScript E}/{\EuScript K}(\tau;{\EuScript I})\langle X_j,X_j^* \mid j \in {\mathbb N}\rangle$ and $F_n(X_1,\dots,X_n) \in {\EuScript E}(\tau;{\EuScript I})\langle X_j,X_j^* \mid j \in {\mathbb N}\rangle$ the non-commutative polynomials
\[
\begin{aligned}
f_n(X_1,\dots,X_n) &= e_n + \sum_{1 \le j \le n} (a_{jn}X_jb_{jn} + c_{jn}X_j^*d_{jn}), \\
F_n(X_1,\dots,X_n) &= E_n + \sum_{1 \le j \le n} (A_{jn}X_jB_{jn} + C_{jn}X_j^*D_{jn}).
\end{aligned}
\]
We shall apply Lemma~3.2 with a sequence $M_k \in {\EuScript E}(\tau;{\EuScript I})$, $k \in {\mathbb N}$ and an increasing function $\varphi: {\mathbb N} \to {\mathbb N}$ such that the set $\{M_k \mid 1 \le k \le \varphi(m)\}$ contains the following operators $E_m$, $A_{jm}$, $B_{jm}$, $C_{jm}$, $D_{jm}$, $Y_{jm}$, $Y_{jm}^*$ where $1 \le j \le m$ and also $F_n(Y_{1m},\dots,Y_{nm})$ and
\[
(F_p(Y_{1m},\dots,Y_{pm}))^*F_q(Y_{1m},\dots,Y_{qm})
\]
with $1 \le n \le m$, $1 \le p \le m$, $1 \le q \le m$. Note that the listed operators won't exhaust $\{M_k \mid 1 \le k \le \varphi(m)\}$, since $\varphi$ being increasing we will have that if $1 \le m' < m$ then $\{M_k \mid 1 \le k \le \varphi(m')\} \subset \{M_k \mid 1 \le k \le \varphi(m)\}$.

Since $p(F_n(Y_{1m},\dots,Y_{nm})) = f_n(y_{1m},\dots,y_{nm})$ if $1 \le n \le m$, we can find $P_k \in {\EuScript P}({\EuScript H})$, $P_k \uparrow I$ so that
\[
|\|F_n(Y_{1m},\dots,Y_{nm})(I-P_m)\| - r_n| < 1/m
\]
if $1 \le n \le m$. Remark that if $1 \le n \le m$ and $N \ge m$ then
\[
\left|\left\|F_n(Y_{1m},\dots,Y_{nm}) \sum_{k \ge N+2} R_k^2\right\| - r_n\right| < 1/m
\]
because $\sum_{k \le N+2} R_k^2 \le I - P_m$ and $I - \sum_{k \ge N+2} R_k^2 \in {\EuScript R}({\EuScript H}) \subset {\EuScript K}({\EuScript H})$ so that
\[
\begin{aligned}
\|F_n(Y_{1m},\dots,Y_{nm})(I-P_m)\| &\ge \left\|F_n(Y_{1m},\dots,Y_{nm}) \sum_{k \ge N+2} R_k^2\right\| \\
&\ge \|f_n(y_{1m},\dots,y_{nm})\|.
\end{aligned}
\]
We can therefore find a sequence $1 < N_1 < N_2 < \dots$ so that $N_m \ge m+2$, $N_{p+1} - N_p \ge 8$ for all $m,p \in {\mathbb N}$ and
\[
\left|\left\|F_n(Y_{1m},\dots,Y_{nm}) \sum_{N_m \le k < N_{m+1}} R_k^2\right\| - r_n\right| < 1/m
\]
if $1 \le n \le m$ and also
\[
\left|\left\|F_n(Y_{1m},\dots,Y_{nm}) \sum_{N_m + 3 \le k < N_{m+1} - 3} R_k^2\right\| - r_n\right| < 1/m.
\]

We will show that if the $\epsilon_m$ are chosen so that $\sum_{m \ge 1} \epsilon_m < \infty$, then the operators
\[
Y_j = \sum_{m \ge j} \left( \sum_{N_m \le k < N_{m+1}} R_kY_{jm}R_k\right)
\]
will satisfy $Y_j \in {\EuScript E}(\tau;{\EuScript I})$, $\|Y_j\| \le 1$ and $p(Y_j) = y_j$ will satisfy $\|f_n(y_1,\dots,y_n)\| = r_n$ for all $n \in {\mathbb N}$.

We will not need to put conditions on the $\epsilon_m$ in order that $\|Y_j\| \le 1$. Indeed, this can be seen as follows. Let $Z: {\EuScript H} \to {\EuScript H} \otimes l^2({\mathbb N})$ be the operator
\[
Zh = \sum_{m \ge j} \sum_{N_m \le k < N_{m+1}} R_kh \otimes e_k
\]
and let $S_j \in {\EuScript B}({\EuScript H} \otimes l^2({\mathbb N}))$ be the operator
\[
S_j \sum_{k \ge 1} h_k \otimes e_k = \sum_{m \ge j} \sum_{N_m \le k < N_{m+1}} Y_{jm}h_k \otimes e_k.
\]
Since $\|Y_{jm}\| < 1$, we have $\|S_j\| \le 1$ and we also have $\|Z\| \le 1$ since 
\[
Z^*Z = \sum_{m \ge j} \sum_{N_m \le k < N_{m+1}} R_k^2 \le \sum_{k \ge 1} R_k^2 = I.
\]
Hence $\|Y_j\| \le 1$ since $Y_j = Z^*S_jZ$.

Our next task will be to show that if $\sum_{m \ge 1} \epsilon_m < \infty$, we will have $|[Y_j,\tau]|_{{\EuScript I}} < \infty$, which together with the boundedness of $Y_j$ we just showed, will give $Y_j \in {\EuScript E}(\tau;{\EuScript I})$.

Since the sum defining $Y_j$ is weakly convergent to $Y_j$, it will be sufficient to show that assuming $\sum_{m \ge 1} \epsilon_m < \infty$ we can insure that
\[
\sum_{m \ge j} \left| \left[ \sum_{N_m \le k < N_{m+1}} R_kY_{jm}R_k,\tau\right]\right|_{{\EuScript I}} < \infty.
\]
Since the $Y_{jm}$ with $1 \le j \le m$ are among the $M_p$ with $1 \le p \le \varphi(m)$ we infer from condition $5^0$ in Lemma~3.2 that $|\|[R_k,Y_{jm}]\|| < \epsilon_k$ if $N_m \le k$ and $1 \le j \le m$. Also by condition $3^0$ of Lemma~3.2, $|\|R_k\|| < 1 + \epsilon_k$. This gives
\[
\begin{aligned}
&\left|\,\left|\left[\sum_{N_m \le k < N_{m+1}} R_kY_{jm}R_k,\tau\right]\right|_{{\EuScript I}} - \left|\left[\sum_{N_m \le k < N_{m+1}} Y_{jm}R_k^2,\tau\right]\right|_{{\EuScript I}}\,\right| \\
&\qquad\le \left|\left\| \sum_{N_m \le k < N_{m+1}} [R_k,Y_{jm}]R_k\right\|\right| \\
&\qquad\le \sum_{N_m \le k < N_{m+1}} |\|[R_k,Y_{jm}]\||\,|\|R_k\|| \\
&\qquad\le \sum_{N_m \le k < N_{m+1}} \epsilon_k(1+\epsilon_k).
\end{aligned}
\]
Hence in order that $Y_j \in {\EuScript E}(\tau;{\EuScript I})$ it will suffice that $\sum_{k \ge 1} \epsilon_k < \infty$ and $\sum_{m \ge j} \left|\left[\sum_{N_m \le k < N_{m+1}} Y_{jm}R_k^2,\tau\right]\right|_{{\EuScript I}} < \infty$. Since $|[R_k,\tau]|_{{\EuScript I}} < \epsilon_k$ if $k \ge 2$ by $3^0$ of Lemma~3.2 and $\|Y_{jm}\| < 1$, we have
\[
\begin{aligned}
\sum_{m \ge j} \left|\left[ \sum_{N_m \le k < N_{m+1}} Y_{jm}R_k^2,\tau\right]\right|_{{\EuScript I}} &\le \sum_{m \ge j} |[Y_{jm},\tau]|_{{\EuScript I}} \cdot \left\| \sum_{N_m\le k < N_{m+1}} R_k^2\right\| \\
&+ \sum_{k \le 2} |[R_k^2,\tau]|_{{\EuScript I}} \le \sum_{m \ge j} \epsilon_m + \sum_{k \ge 2} 2\epsilon_k < \infty
\end{aligned}
\]
under the assumption that $\sum_{k \ge 1} \epsilon_k < \infty$. Hence under this condition on the $\epsilon_m$ we have $Y_j \in {\EuScript E}(\tau;{\EuScript I})$.

Finally, we turn to showing that assuming $\sum_{m \ge 1} \epsilon_m < \infty$, we will have $\|f_n(y_1,\dots,y_n)\| = r_n$ for all $n \in {\mathbb N}$, where $y_j = p(Y_j)$. Clearly $f_n(y_1,\dots,y_n) = p(F_n(Y_1,\dots,Y_n))$. Note also that the relations we're aiming at being about norms in the Calkin algebra, we will no longer have to deal with $|\|\cdot\||$-norms and the ideal ${\EuScript I}$ for this matter.

We begin by showing that we can arrange that the difference between $F_n(Y_1,\dots,Y_n)$ and
\[
\sum_{m \ge n} \sum_{N_m \le k < N_{m+1}} R_kF_n(Y_{1m},\dots,Y_{nm})R_k
\]
is a compact operator. Since
\[
F_n(Y_1,\dots,Y_n) = E_n + \sum_{1 \le j \le n} (A_{jn}Y_jB_{jn}+C_{jn}Y_j^*D_{jn})
\]
it will suffice to prove the assertion in 3 cases, when $F_n(Y_1,\dots,Y_n)$ equals $E_n$, $A_{jn}Y_jB_{jn}$, $C_{jn}Y_{jn}^*D_{jn}$ where $1 \le j \le n$.

In the first case we have
\[
\sum_{k \ge N_n} R_kE_nR_k - E_n = -\sum_{1 \le k < N_n} R_kE_nR_k + \sum_{k \ge 1} (R_kE_nR_k - E_nR_k^2).
\]
The first sum being finite rank, we need that the second sum be compact. If $k \ge n$, $\|[R_k,E_n]\| < \epsilon_k$ since $E_n$ is among the $M_p$ with $p \le \varphi(n) \le \varphi(k)$ and condition $5^0$ of Lemma~3.2 holds. Thus, $\|R_kE_nR_k - E_nR_k^2\| < \epsilon_k$ implies that the difference we consider will be compact if $\sum_{k \ge 1} \epsilon_k < \infty$.

In case $F_n$ is $A_{jn}Y_jB_{jn}$, where $1 \le j \le n$, we must insure compactness of
\[
\begin{aligned}
&\sum_{m \ge n} \sum_{N_m \le k < N_{m+1}} R_kA_{jn}Y_{jm}B_{jn}R_k - \sum_{m \ge n} A_{jn}\left( \sum_{N_m \le k < N_{m+1}} R_kY_{jm}R_k\right)B_{jn} \\
&\qquad = \sum_{m \ge n} \sum_{N_m \le k < N_{m+1}} ([R_k,A_{jn}]Y_{jm}B_{jn}R_k + A_{jn}R_kY_{jm}[B_{jn},R_k]).
\end{aligned}
\]
The last sum being a sum of finite rank operators it will suffice to have convergence of the sum of their norms. Since the $A_{jn}$ and $B_{jn}$ are among the $M_p$ with $p \le \varphi(n) \le \varphi(k)$ we have that the norms of the commutators are majorized by $\epsilon_k$ in view of $5^0$ in Lemma~3.2 and hence the sum of norms is majorized by
\[
K \sum_{k \ge 1} \epsilon_k
\]
where $K$ is a bound for $\|A_{jn}\|$ and $\|B_{jn}\|$. Thus again it will suffice that $\sum_{k \ge 1} \epsilon_k < \infty$.

The third situation when we consider $C_{jn}Y_j^*D_{jm}$ is entirely analogous to that of $A_{jn}Y_jB_j$, since we treated $Y_{jn}$ and $Y_{jn}^*$ symmetrically in our assumptions about $\varphi$. Again, summability of the $\epsilon_m$ will suffice.

We need then to show that if $\sum_{m \ge 1} \epsilon_m < \infty$, we will also have that the essential norm of
\[
\Omega_n = \sum_{m \ge n} \sum_{N_m \le k < N_{m+1}} R_kF_n(Y_{1m},\dots,Y_{nm})R_k
\]
will be $r_n$.

Using again the operator
\[
Z: {\EuScript H} \to {\EuScript H} \otimes l^2({\mathbb N}),\ Zh = \sum_{k \ge N_n} R_kh \otimes e_k
\]
we have $\|Z\| \le 1$ and $Z^*\Gamma_{nt}Z - \Omega_n \in {\EuScript R}({\EuScript H})$, where for $t \ge n$ we define on ${\EuScript H} \otimes l^2({\mathbb N})$ operators $\Gamma_{nt}$ by
\[
\Gamma_{nt}\sum_{k \ge 1} h_k\otimes e_k = \sum_{m \ge t} \sum_{N_m \le k < N_{m+1}} F_n(Y_{1m},\dots,Y_{nm})(I-P_m) h_k \otimes e_k.
\]
Since $\|\Gamma_{nt}\| = \sup_{m \ge t} \|F_n(Y_{1m},\dots,Y_{nm})(I-P_m)\|$ we have $|\|\Gamma_{nt}\| - r_n| < t^{-1}$ and hence $\lim_{t \to +\infty} \|\Gamma_{nt}\| = r_n$. This gives $\|p(\Omega_n)\| \le \lim_{t \to +\infty}\|\Gamma_{nt}\| = r_n$ and hence we are left with the opposite inequality $\|p(\Omega_n)\| \le r_n$.

We will again use a compact perturbation and pass from $\Omega_n$ to another operator
\[
\Xi_n = \sum_{m \ge n} F_n(Y_{1m},\dots,Y_{nm}) \sum_{N_m \le k < N_{m+1}} R_k^2.
\]
Indeed we have
\[
\Xi_n - \Omega_n = \sum_{m \ge n} \sum_{N_m \le k < N_{m+1}} [F_n(Y_{1m},\dots,Y_{nm}),R_k]R_k
\]
and
\[
\|[R_n(Y_{1m},\dots,Y_{nm}),R_k]R_k\| < \epsilon_k.
\]
Again compactness will follow if $\sum_{k \ge 1} \epsilon_k < \infty$.

Recall now that we had chosen the $N_m$ so that
\[
\left|\left\| F_n(Y_{1m},\dots,Y_{nm}) \sum_{N_m \le k < N_{m+1}} R_k^2\right\| - r_n\right| < 1/m
\]
and also
\[
\left|\left\| F_n(Y_{1m},\dots,Y_{nm}) \sum_{N_m+3 \le k < N_{m+1}-3} R_k^2\right\| - r_n\right| < 1/m.
\]
Since by Lemma~3.2 we have that the $R_k$ are finite rank positive contractions, commute and satisfy $|k-l| \ge 2 \Rightarrow R_kR_l = 0$ it is easily seen that if $\Delta_m$ is the projection onto the range of $\sum_{N_m+3 \le k < N_{m+1}-3} R_k^2$ we will have $R_s\Delta_m = 0$ if $s < N_m$ or $s \ge N_{m+1}$ and hence
\[
\Delta_p\Delta_q = 0 \text{ if } p \ne q
\]
and if $n \le m$ we have
\[
\begin{aligned}
\left\| F_n(Y_{1m},\dots,Y_{nm})\sum_{N_m+3 \le k < N_{m+1}-3} R_k^2\right\| &\le \|\Xi_n\Delta_m]\| \\
&\le \left\| F_n(Y_{1m},\dots,Y_{nm}) \sum_{N_m \le k < N_{m+1}} R_k^2\right\|
\end{aligned}
\]
so that
\[
|\|\Xi_n\Delta_m\| - r_n| < 1/m.
\]
This implies
\[
\|p(\Xi_n)\| \ge \limsup_{m \to +\infty} \|\Xi_n\Delta_m\| = r_n.
\]
\qed

\section{Multipliers and duality}
\label{sec4}

In this section we will sometimes also deal with normed ideals in which the finite rank operators are not dense, which occurs when the norming function $\Phi$ is not mononorming (see the preliminaries and \cite{9} or \cite{12}). We begin with a basic lemma.

\bigskip
\noindent
{\bf 4.1. Lemma.} {\em Let $\Phi$ be a norming function and let $({\EuScript I},|\ |_{{\EuScript I}}) = ({\textgoth S}_{\Phi},|\ |_{\Phi})$ so that $({\EuScript I}^{(0)},|\ |_{{\EuScript I}}) = ({\textgoth S}_{\Phi}^{(0)},|\ |_{\Phi})$ is the closure of ${\EuScript R}({\EuScript H})$ in ${\EuScript I}$. Assume $k_{{\EuScript I}}(\tau) = 0$. Then ${\EuScript K}(\tau;{\EuScript I}^{(0)})$ is a closed two-sided ideal in ${\EuScript E}(\tau;{\EuScript I})$ and the norm in ${\EuScript E}(\tau;{\EuScript I})$ extends the norm in ${\EuScript K}(\tau;{\EuScript I}^{(0)})$. Moreover, the unit ball of $({\EuScript E}(\tau;{\EuScript I}),|\|\cdot\||)$ is weakly compact.
}

\bigskip
\noindent
{\bf {\em Proof.}} It is clear that the norm of ${\EuScript E}(\tau;{\EuScript I})$ extends the norm of ${\EuScript K}(\tau;{\EuScript I}^{(0)})$ and that ${\EuScript K}(\tau;{\EuScript I}^{(0)})$ is a closed subalgebra of ${\EuScript E}(\tau;{\EuScript I})$. By Corollary~2.2 ${\EuScript R}({\EuScript H})$ is dense in ${\EuScript K}(\tau;{\EuScript I}^{(0)})$ and hence ${\EuScript K}(\tau;{\EuScript I}^{(0)})$ is the closure in ${\EuScript E}(\tau;{\EuScript I})$ of the two-sided ideal ${\EuScript R}({\EuScript H})$, which implies that also ${\EuScript K}(\tau;{\EuScript I}^{(0)})$ is a two-sided ideal in ${\EuScript E}(\tau;{\EuScript I})$. If $X$ is the weak limit of the net $(X_{\alpha})_{\alpha \in I}$ in the unit ball of ${\EuScript E}(\tau;{\EuScript I})$ then by the weak compactness of the unit balls of ${\EuScript B}({\EuScript H})$ and of ${\textgoth S}_{\Phi}$ (see \cite{9}) we have $\|X\| \le 1$ and $|[X,\tau]|_{\Phi} \le 1$ so that $X \in {\EuScript E}(\tau;{\EuScript I})$. Since ${\EuScript H}$ is separable we may replace $(X_{\alpha})_{\alpha \in I}$ by a subsequence and use the semicontinuity properties of $\|\ \|$ and $|\ |_{\Phi}$ under weak convergence to get that $|\|X\|| \le 1$. Thus the unit ball of ${\EuScript E}(\tau;{\EuScript I})$ is a closed subset of the unit ball of ${\EuScript B}({\EuScript H})$ and hence weakly compact.\qed

\bigskip
We pass now to bounded multipliers ${\EuScript M}({\EuScript K}(\tau;{\EuScript I}^{(0)}))$ that is double centralizer pairs $(T',T'')$ of bounded linear maps ${\EuScript K}(\tau;{\EuScript I}^{(0)}) \to {\EuScript K}(\tau;{\EuScript I}^{(0)})$ so that $T'(x)y = xT''(y)$ (\cite{10}).

\bigskip
\noindent
{\bf 4.2. Proposition.} {\em Assume $k_{{\EuScript I}}(\tau) = 0$, where ${\EuScript I} = {\textgoth S}_{\Phi}$ and ${\EuScript I}^{(0)} = {\textgoth S}_{\Phi}^{(0)}$. We have ${\EuScript M}({\EuScript K}(\tau;{\EuScript I}^{(0)})) = {\EuScript E}(\tau;{\EuScript I})$, that is, if $(T',T'') \in {\EuScript M}({\EuScript K}(\tau;{\EuScript I}^{(0)}))$ then there is a unique $T \in {\EuScript E}(\tau;{\EuScript I})$ so that $T'(x) = xT$ and $T''(x) = Tx$.
}

\bigskip
\noindent
{\bf {\em Proof.}} By Corollary~2.4 there is a sequence $A_s \in {\EuScript R}_1^+({\EuScript H})$ so that $\|A_s\| = 1$, $s > t \Rightarrow A_sA_t = A_t$ and $A_s \uparrow I$, $|\|A_s\|| \to 1$, $|\|(I-A_s)K\|| \to 0$ if $s \to \infty$ and $K \in {\EuScript K}(\tau;{\EuScript I}^{(0)})$.

Assume $(T',T'') \in {\EuScript M}({\EuScript K}(\tau;{\EuScript I}^{(0)}))$ and let $K_s = T'(A_s)A_s = A_sT''(A_s)$. Clearly $\sup_{s \in {\mathbb N}} |\|K_s\|| < \infty$ the multiplier being bounded. Remark also that $s > t \Rightarrow A_tK_sA_t = A_tT'(A_s)A_sA_t = A_tT'(A_s)A_t = A_tA_sT''(A_t) = A_tT''(A_t) = K_t$. Hence if $T$ is the weak limit of a subsequence of the $K_s$, we have $A_tTA_t = K_t$ for all $t$ and hence $T$ does not depend on the subsequence, that is $T = w - \lim_{s \to \infty} K_s$ and also $T \in {\EuScript E}(\tau;{\EuScript I})$ since the unit ball of ${\EuScript E}(\tau;{\EuScript I})$ is weakly closed.

On the other hand if $K \in {\EuScript K}(\tau;{\EuScript I}^{(0)})$ then $|\|A_sK - K\|| \to 0$ as $s \to \infty$ and also $|\|KA_s - K\|| \to 0$ as $s \to \infty$ (replace $K$ by $K^*$). We have
\[
\begin{aligned}
T'(K)A_t &= \lim_{s \to \infty} T'(KA_s)A_t \\
&= \lim_{s \to \infty} KA_sT''(A_t) \\
&= \lim_{s \to \infty} KT'(A_s)A_t \\
&= \lim_{s \to \infty} KT'(A_s)A_sA_t = KTA_t
\end{aligned}
\]
and since this holds for all $t \in {\mathbb N}$ we have $T'(K) = KT$. This then gives $T'(A_t)K = A_tTK = A_tT''(K)$ and hence $T''(K) = TK$.

Uniqueness of $T$ follows from ${\EuScript K}(\tau;{\EuScript I}^{(0)}) \supset {\EuScript R}({\EuScript H})$. The converse, that $T \in {\EuScript E}(\tau;{\EuScript I})$ gives rise to a multiplier, is a consequence of Lemma~4.1.\qed

\bigskip
We pass now to duality. Recall from the theory of normed ideals (\cite{9}, \cite{12}) that given norming function $\Phi$ there is a conjugate norming function $\Phi^*$ so that the dual of the Banach space $({\textgoth S}_{\Phi}^{(0)},|\ |_{\Phi})$ is $({\textgoth S}_{\Phi^*},|\ |_{\Phi^*})$ under the duality $(X,Y) \to \text{\rm Tr}(XY)$ for $(X,Y) \in {\textgoth S}_{\Phi}^{(0)} \times {\textgoth S}_{\Phi^*}$ (we leave out of the discussion the case of ${\textgoth S}_{\Phi}^{(0)} = {\EuScript C}_1$, where the dual is ${\EuScript B}({\EuScript H})$).

\bigskip
\noindent
{\bf 4.3. Proposition.} {\em Let $\Phi$ be a norming function so that $k_{\Phi}(\tau) = 0$, let $\Phi^*$ be its conjugate and assume ${\textgoth S}_{\Phi}^{(0)} \ne {\EuScript C}_1$. Then the dual of ${\EuScript K}(\tau;{\textgoth S}_{\Phi}^{(0)})$ can be identified isometrically with $({\EuScript C}_1 \times ({\textgoth S}_{\Phi^*})^n)/{\EuScript N}$ where
\[
\begin{split}
{\EuScript N} = \left\{ \left( \sum_{1 \le j \le n} [T_j,y_j],(y_j)_{1 \le j \le n}\right) \in {\EuScript C}_1 \times ({\textgoth S}_{\Phi^*})^n \mid \right. \\
\left. (y_j)_{1 \le j \le n} \in ({\textgoth S}_{\Phi^*})^n \text{ with } \sum_{1 \le j \le n} [T_j,y_j] \in {\EuScript C}_1\right\}
\end{split}
\]
and the duality map ${\EuScript K}(\tau;{\textgoth S}_{\Phi}^{(0)}) \times ({\EuScript C}_1 \times ({\textgoth S}_{\Phi^*})^n) \to {\mathbb C}$ is 
\[
(K,(x,(y_j)_{1 \le j \le n})) \to \text{\rm Tr}(Kx + \sum_{1 \le j \le n} [T_j,K]y_j)
\]
and the norm on $({\EuScript C}_1 \times ({\textgoth S}_{\Phi^*})^n)$ is 
\[
\|(x,(y_j)_{1 \le j \le n})\| = \max\left( |x|_1,\sum_{1 \le j \le n} |y_j|_{\Phi^*}\right).
\]
}

\bigskip
\noindent
{\bf {\em Proof.}} Since $K \to K \oplus [\tau,K]$ identifies ${\EuScript K}(\tau;{\textgoth S}_{\Phi}^{(0)})$ isometrically with a closed subspace of ${\EuScript K}({\EuScript H}) \oplus ({\textgoth S}_{\Phi}^{(0)})^n$ with the norm $\|K \oplus (H_j)_{1 \le j \le n}\| = \|K\| + \max_{1 \le j \le n} |H_j|_{\Phi}$, the dual of which is ${\EuScript C}_1 \times ({\textgoth S}_{\Phi^*})^n$, the proof boils down to showing that ${\EuScript N}$ is the annihilator of
\[
\{K \oplus [\tau,K] \in {\EuScript K}({\EuScript H}) \oplus ({\textgoth S}_{\Phi}^{(0)})^n \mid K \in {\EuScript K}(\tau;{\textgoth S}_{\Phi}^{(0)})\}.
\]
Since ${\EuScript R}({\EuScript H})$ is dense in ${\EuScript K}(\tau;{\textgoth S}_{\Phi}^{(0)})$ by Corollary~2.2, it will suffice to show that ${\EuScript N}$ is the annihilator of
\[
\{R \oplus [\tau,R] \in {\EuScript K}({\EuScript H}) \oplus ({\textgoth S}_{\Phi}^{(0)})^n \mid R \in {\EuScript R}({\EuScript H})\}.
\]
If $R \in {\EuScript R}({\EuScript H})$ and $(x,(y_j)_{1 \le j \le n}) \in {\EuScript N}$ we have
\[
\begin{aligned}
\text{\rm Tr}\left(Rx + \sum_{1 \le j \le n} [T_j,R]y_j\right) &= \text{\rm Tr}\left( R \sum_{1 \le j \le n} [T_j,y_j] + \sum_{1 \le j \le n} [T_j,R]y_j\right) \\
&= \text{\rm Tr}\left( \sum_{1 \le j \le n} [T_j,Ry_j]\right) = 0.
\end{aligned}
\]
Conversely, if $(x,(y_j)_{1 \le j \le n}) \in {\EuScript C}_1 \times ({\textgoth S}_{\Phi^*})^n$ is such that
\[
\text{\rm Tr}\left(Rx + \sum_{1 \le j \le n} [T_j,R]y_j\right) = 0 \text{ for all } R \in {\EuScript R}({\EuScript H})
\]
then
\[
\text{\rm Tr}\left(R\left( x - \sum_{1 \le j \le n} [T_j,y_j]\right)\right) = 0 \text{ for all } R \in {\EuScript R}({\EuScript H})
\]
and hence $x = \sum_{1 \le j \le n} [T_j,y_j]$ that is $(x,(y_j)_{1 \le j \le n}) \in {\EuScript N}$.\qed

\bigskip
\noindent
{\bf 4.4. Lemma.} {\em Under the assumptions of Proposition~$4.3$ and the additional assumption that $\Phi^*$ is mononorming
\[
\left\{\left( \sum_{1 \le j \le n} [T_j,R_j],(R_j)_{1 \le j \le n}\right) \in {\EuScript C}_1 \times ({\textgoth S}_{\Phi^*})^n \mid (R_j)_{1 \le j \le n} \in ({\EuScript R}({\EuScript H}))^n\right\}
\]
is dense in ${\EuScript N}$.
}

\bigskip
\noindent
{\bf {\em Proof.}} Let $(x,(y_j)_{1 \le j \le n}) \in {\EuScript N}$, that is $(y_j)_{1 \le j \le n} \in ({\textgoth S}_{\Phi^*})^n$ is such that $x = \sum_{1 \le j \le n} [T_j,y_j] \in {\EuScript C}_1$. Let $A_s \in {\EuScript R}_1^+({\EuScript H})$, $s \in {\mathbb N}$ be such that $A_s \uparrow I$ and $|[\tau,A_s]|_{\Phi} \to 0$ as $s \to \infty$. Since ${\textgoth S}_{\Phi^*} = {\textgoth S}_{\Phi^*}^{(0)}$, $\Phi^*$ being mononorming, we have $|y_jA_s - y_j|_{\Phi^*} \to 0$ as $s \to \infty$, $1 \le j \le n$. Moreover, we have
\[
\begin{aligned}
&\left| \sum_{1 \le j \le n} [T_j,y_jA_s] - \sum_{1 \le j \le n} [T_j,y_j]\right|_1 \\
&\qquad = \left|
 \left( \sum_{1 \le j \le n} [T_j,y_j]\right)A_s + \sum_{1 \le j \le n} y_j[T_j,A_s] - \sum_{1 \le j \le n} [T_j,y_j]\right|_1 \\
&\qquad\le \left| \left( \sum_{1 \le j \le n} [T_j,y_j]\right)(I-A_s)\right|_1 + \sum_{1 \le j \le n} |y_j|_{\Phi^*}|[T_j,A_s]|_{\Phi} \to 0
\end{aligned}
\]
as $s \to \infty$.\qed

\bigskip
\noindent
{\bf 4.5. Proposition.} {\em Assume $k_{\Phi}(\tau) = 0$ and that both norming functions $\Phi$ and $\Phi^*$ are mononorming. Then with the notation of Proposition~$4.3$, the dual of $({\EuScript C}_1 \times ({\textgoth S}_{\Phi^*})^n)/{\EuScript N}$ identifies with ${\EuScript E}(\tau;{\textgoth S}_{\Phi})$ via the duality map
\[
\left(X,(x,(y_j)_{1 \le j \le n})\right) \to \text{\rm Tr}\left(Xx + \sum_{1 \le j \le n} [T_j,X]y_j\right).
\]
In particular, ${\EuScript E}(\tau;{\textgoth S}_{\Phi})$ identifies with the bidual of ${\EuScript K}(\tau;{\textgoth S}_{\Phi})$ (note that under our assumptions ${\textgoth S}_{\Phi}^{(0)} = {\textgoth S}_{\Phi}$).
}

\bigskip
\noindent
{\bf {\em Proof.}} The dual of $({\EuScript C}_1 \times ({\textgoth S}_{\Phi^*})^n)/{\EuScript N}$ is the orthogonal of ${\EuScript N}$ in the dual of ${\EuScript C}_1 \times ({\textgoth S}_{\Phi^*})^n$. The mononorming assumption on $\Phi$ and $\Phi^*$ implies ${\textgoth S}_{\Phi}$ and ${\textgoth S}_{\Phi^*}$ are each others dual and are reflexive and separable (see Thm.~12.2 in Ch.~III of \cite{9}). The dual of $({\EuScript C}_1 \times ({\textgoth S}_{\Phi^*})^n)$ is ${\EuScript B}({\EuScript H}) \oplus ({\textgoth S}_{\Phi})^n$ (the usual duality based on the trace). Since Lemma~4.4 provides a dense subset of ${\EuScript N}$, it suffices to show that $\{X \oplus [\tau,X] \in {\EuScript B}({\EuScript H}) \oplus ({\textgoth S}_{\Phi})^n \mid X \in {\EuScript E}(\tau;{\textgoth S}_{\Phi})\}$ is the orthogonal in ${\EuScript B}({\EuScript H}) \oplus ({\textgoth S}_{\Phi})^n$ of the set
\[
\left\{\left( \sum_{1 \le j \le n} [T_j,R_j],(R_j)_{1\le j \le n}\right) \in {\EuScript C}_1 \times ({\textgoth S}_{\Phi^*})^n \mid (R_j)_{1 \le j \le n} \in ({\EuScript R}({\EuScript H}))^n\right\}.
\]
Indeed, if $X \oplus (H_j)_{1 \le j \le n} \in {\EuScript B}({\EuScript H}) \oplus ({\textgoth S}_{\Phi})^n$ is such that \[
\text{\rm Tr}\left(X \sum_{1 \le j \le n}[T_j,R_j] + \sum_{1 \le j \le n} R_jH_j\right) = 0
\]
for all $(R_j)_{1 \le j \le n} \in ({\EuScript R}({\EuScript H}))^n$ then $\text{\rm Tr}\left( \sum_{1 \le j \le n} ([X,T_j] + H_j)R_j\right) = 0$ for all $(R_j)_{1\le j \le n} \in ({\EuScript R}({\EuScript H}))^n$, which implies $H_j = [T_j,X]$, $1 \le j \le n$. Thus we have proved that $[T_j,X] \in {\textgoth S}_{\Phi}$ and hence also $X \in {\EuScript E}(\tau;{\textgoth S}_{\Phi})$. Clearly, also if $X \in {\EuScript E}(\tau;{\textgoth S}_{\Phi})$ and $(R_j)_{1 \le j \le n} \in ({\EuScript R}({\EuScript H}))^n$ we have
\[
\text{\rm Tr}\left( X \sum_{1 \le j \le n} [T_j,R_j] + \sum_{1 \le j \le n} [T_j,X]R_j\right) = \text{\rm Tr}\left( \sum_{1 \le j \le n} [T_j,XR_j]\right) = 0.
\]
\qed

\section{Concluding remarks}
\label{sec5}

One may wonder what happens with several of the results of this paper if some of the conditions ${\EuScript I} = {\EuScript I}^{(0)}$ and $k_{{\EuScript I}}(\tau) = 0$ are relaxed or removed. For instance, it is natural to ask whether the density of ${\EuScript R}({\EuScript H})$ in ${\EuScript K}(\tau;{\EuScript I})$ when ${\EuScript I} \ne {\EuScript I}^{(0)}$ and $k_{{\EuScript I}}(\tau) = 0$ is possible. Other questions, such as studying ${\EuScript E}/{\EuScript K}(\tau;{\EuScript I})$ when $k_{{\EuScript I}}(\tau) \ne 0$ or ${\EuScript I} \ne {\EuScript I}^{(0)}$ may require quite different methods in case ${\EuScript E}/{\EuScript K}(\tau;{\EuScript I})$ is no longer a $C^*$-algebra.

\end{document}